\documentclass[a4paper, 12pt]{article}

\usepackage{graphics,color}
\def\spc{\color[rgb]{1,0.2,0.2}}
\def\nc{\normalcolor}

\usepackage{amsmath, amssymb}
\usepackage{verbatim}
\usepackage[british]{babel}
\usepackage{graphicx}

\usepackage{color}

\newtheorem {definition}{Definition}
\newtheorem {lemma}{Lemma}

\newtheorem {prop}{Proposition}
\newtheorem {theorem}{Theorem}
\newtheorem{rem}{Remark}

\newcommand{\n}{{\rm n}}

\newcommand{\R}{{\mathbb{R}}}
\newcommand{\N}{{\mathbb{N}}}
\renewcommand{\P}{\mathsf{P}}

\newcommand{\Pb}{\mathbb{P}}
\newcommand{\Qb}{\mathbb{Q}}

\newcommand{\C}{\mathsf{C}}
\newcommand{\KInC}{\mathsf{C_{{\rm in}}}}
\newcommand{\KOutC}{\mathsf{C_{{\rm out}}}}

\newcommand{\Put}{\mathsf{Put}}

\newcommand{\BSCall}{\mathsf{BSC}}
\newcommand{\BSPut}{\mathsf{BSP}}

\newcommand{\Z}{\mathbb{Z}}
\newcommand{\eps}{\varepsilon}
\newcommand{\E}{\mathsf{E}}

\newcommand{\Xn}{X^{\left(n\right)}}
\newcommand{\Sn}{S^{\left(n\right)}}
\newcommand{\Wp}{W^{\left(p\right)}}

\begin{document}

\title{Density of Skew Brownian motion and  its  functionals \\ with  application in  finance }

\author{
Alexander Gairat\footnote{ Gazprombank, Moscow.
 Email address: Alexander.Gairat@gazprombank.ru
}\, and 
Vadim Shcherbakov\footnote{
 Department of Mathematics, Royal Holloway,  University of London.
 Email address: vadim.shcherbakov@rhul.ac.uk
}
}

\date{}

\maketitle

\begin{abstract}
{\small 
We  derive the joint density of a Skew Brownian motion, its last visit to the origin,  its  local and occupation times.
 The result  allows to obtain explicit analytical formulas  for  pricing  European options 
  under both  a two valued local volatility model and a displaced diffusion model
with constrained volatility.
}

\medskip 

{\small {\bf Key words}: Skew Brownian motion, local volatility model, displaced diffusion, 
 local time, occupation time, simple random walk, option pricing} 
\end{abstract}

\section{Introduction}
\label{intro}
 
A Skew Brownian motion (SBM) with parameter $p$  is a Markov process that  evolves  as 
 a standard Brownian motion  reflected at the origin so that the next excursion is chosen to be positive 
 with probability $p$.   SBM   was introduced in  \cite{Ito} and has been  studied extensively   in probability since  then.
The process naturally appears  in diverse  applications, e.g.   \cite{Appu} and  \cite{Lejay},  and, in particular,
  in finance applications, e.g.  \cite{Decamps0}, \cite{Decamps1}, \cite{Decamps2} and \cite{Rossello}.
In this paper we derive the joint distribution of SBM and some of its functionals and apply this distribution 
 to   derivative pricing under both  a local volatility model with discontinuity 
and a displaced diffusion model  with constrained volatility.

Let $(\Omega, {\cal F}, \P)$ be a probability space and  let  $\{W_t, {\cal F}_t,\, t\geq 0\}$ be a  standard Brownian motion (BM)
with its natural  filtration.  As usual, denote by $\R$ and $\R_{+}$ sets of all real and all non-negative real numbers respectively.
A local volatility model (LVM) for the underlying price $S_t$ is given by the following equation 
\begin{equation}
\label{S0}
dS_t=\mu(t)S_tdt+ \sigma(t, S_t)S_tdW_t, 
\end{equation}
where $\mu(t)\in \R$ and $\sigma(t, S_t)\in \R_{+}$.
LVM is a natural extension of the famous Black-Scholes model. The latter 
is a particular case of~(\ref{S0})  where   both drift $\mu$  and volatility $\sigma$ are constant.
LVM  is  actively used in practice because it can be easily calibrated to the market. 
Furthermore,  by Gyongy's lemma (\cite{Gyongy}) a wider class of stochastic volatility models can be reduced to LVM.

A  number of approximations to LVM have been   developed
 for  both calibration purposes and qualitative  analysis  (\cite{Guyon}). 
We apply our probabilistic results  primarily 
to  a  particular case of  LVM  that can be used as a benchmark model for analyzing quality of such approximations.
Namely, we consider  a  driftless LVM  with  a two-valued  volatility (two-valued LVM) 
\begin{equation}
\label{sigma0}
\sigma(t, S)=\sigma_11_{\{S\geq S^*\}}+\sigma_21_{\{ S<S^*\}},
\end{equation}
where $\sigma_i>0, i=1,2,\, S^*>0$, and $1_A$ is used to denote the indicator function of set $A$.
Without loss of generality we assume that $S^*=1$ in what follows. 

In Section~\ref{relation} we show  that if  $S_t$ follows  the two-valued LVM  then  rescaled process $X_t=\log(S_t)/\sigma(S_t)$ 
is a solution of a stochastic differential equation (SDE) of the following type
\begin{equation}
\label{X0}
X_t=X_0+ \int\limits_0^tm(X_s)ds+(2p-1)L_t^{\left(0\right)}(X)+W_t,
\end{equation}
where  $L_t^{\left(0\right)}(X)$ 
is the  local time of process $X_t$  at zero, $p\in (0, 1)$, 
\begin{equation}
\label{m}
m(x)=m_11_{\{x\geq 0\}}+m_21_{\{x<0\}},\, m_1, m_2\in \R,
\end{equation}
and  both $p$ and pair $(m_1, m_2)$ are uniquely determined by $\sigma_1$ and $\sigma_2$ (Lemma~\ref{T0}). Notice that  
SDE~(\ref{X0}) belongs to the following  class of SDE with local time
\begin{equation}
\label{X2}
dX_t=b(X_t)dt +\sigma(X_t)dW_t+\int\limits_{\R}\nu(dx)dL_t^{\left(x\right)}(X),
\end{equation}
where $\nu$ is a finite  signed measure with atoms at the points, where both $b$ and  $\sigma$ can be discontinuous, and
$L_t^{\left(x\right)}(X)$ is the  local time of process $X$ at $x$.
It is known  that  SDE~(\ref{X2})  has a unique  strong solution under certain general conditions  which 
 are satisfied   in the case of equation~(\ref{X0})  (e.g. \cite{LeGall}, \cite{Lejay} and references therein).
In particular, if $m\equiv 0$ then a unique strong solution  of equation (\ref{X0}) is  a SBM with parameter $p$
which we are going to denote  by $W^{\left(p\right)}_t$ from now on.
If $m_1=m_2=m$, then equation (\ref{X0}) takes the following form 
\begin{equation}
\label{SBMD}
X_t=X_0+mt +(2p-1)L_t^{\left(0\right)}+W_t.
\end{equation}
A diffusion process defined by  equation (\ref{SBMD})
appears, for instance,  in a  study of dispersion across an interface  in   \cite{Appu}
and  is named  there  as a SBM with parameter $p$ and drift $m$.  
By analogy, we refer to the solution of equation (\ref{X0}) as a  SBM  with  a two-valued drift.
A SBM  with  two-valued drift (\ref{m})  is  reflected at the origin in the same way as a driftless  $W^{\left(p\right)}_t$ 
and evolves  as a  BM with drift $m_1$ when  it is  above zero  and, respectively,  with  drift $m_2$ when it is below zero.

In general, option prices under  LVM  are calculated by solving numerically  the corresponding partial differential equations,
though some  semi-analytical results are also known. 
For example, semi-analytical results have been obtained in~\cite{Decamps1} for LVM where 
 $\sigma(t, S)=\sigma(S)$ is continuous at all but one point.  
Another example  is provided by~\cite{Sepp}, where   semi-analytical results  have been obtained for 
 LVM with a so-called tiled local volatility. 

We  derive explicit formulas for  the joint density  of  $W^{\left(p\right)}_t$, its last visit to the origin,  local and occupation times
in both the driftless  and the two-valued drift cases.
The joint density  is then  applied    to option pricing    under both   LVM  with volatility~(\ref{sigma0})
and  a displaced diffusion model with constrained volatility. 
The second model   is defined in Section~\ref{displace} and  is of a particular interest  in practical applications. 
It should be noted that both models belong to the  more  general class of LVM considered in~\cite{Decamps1}.
It turns out that European option prices in both cases can be expressed analytically in terms of both  a standard univariate 
 normal distribution and  a bivariate normal distribution. 

Joint distributions of SBM and its functionals are  of interest on their  own right. 
For example, the joint density   of SBM with a constant drift, its local and occupation times was obtained   in~\cite{Appu}.
The result of~\cite{Appu}  generalizes the classic  result of~\cite{Karatzas1}, where the same trivariate density  was obtained 
for the standard BM.
In~\cite{Appu}   the technique of  \cite{Ito} was modified  to 
obtain a Feynman--Kac formula for SBM and this   allowed them   to adopt  the method 
of~\cite{Karatzas1}. In turn, the  method of~\cite{Karatzas1} is  based on the computation of the Laplace transform 
of the joint density. 
In contrast,  we use  a  discrete approximation of SBM by a random walk
 and a key step of our approach consists in 
combining  an intuitively clear path  decomposition for  the  discrete process 
with some well known properties of the symmetric simple random walk.  This  allows us    to derive  analytically  tractable
expressions for the joint density  of discrete analogues of quantities of interest 
and to compute the limit density.    
 
 A discrete approximation is a well known  method  for obtaining
joint distribution of both BM and SBM and their functionals  (e.g. \cite{Lulko} or  \cite{Tak}).
 We were inspired by the use of this method  in  \cite{Bill} 
 for  computation   of  the joint distribution of the  standard BM, its both the  occupation time and the  last visit to the origin.

The paper is organized as follows. We formulate the results on the joint distribution 
of SBM and its functionals in Section~\ref{dens}. In  Section~\ref{finance} we describe 
the relationship between  LVM  with the two-valued volatility   and SBM with the two-valued drift.
Theorem~\ref{T2} in Section~\ref{PricingTh} is an example of an option pricing theorem under the two-valued   LVM. Proofs are given in Section~\ref{proofs}.
In Section~\ref{simple} we also derive a simple closed form approximation for option prices based on the Black-Scholes formula. Effectiveness of the approximate result is tested in comparison with the exact result provide by Theorem~\ref{T2} and another approximation obtained in~\cite{Sepp}. 
Finally, we  discuss in Section~\ref{displace} how our results can be applied to derivative pricing under the 
displaced  diffusion model with constrained volatility.

\section{Density of Skew Brownian motion,  its last visit to the origin, occupation and local times}
\label{dens}

Given a continuous semimartingale   $X_t,\, t\in [0,T]$, define 
the following quantities
\begin{align}
\label{tau0}
\tau_0&= \min \left\{t:X_t=0\right\}\\
\label{tau}
\tau&=\max\{t\in (0,T]: X_t=0\},\\
\label{V}
V&=\int\limits_{\tau_0}^{\tau}1_{\{X_t\geq 0\}}dt.
\end{align}
Also,  let
$L_t^{\left(x\right)}(X)$ be the symmetric  local time  of $X_t$   at point $x$.
In particular,  if  $X_t=W^{\left(p\right)}_t$,  or $X_t$ is the  solution of equation~(\ref{X0}) then
\begin{equation}
\label{L}
L_t^{\left(x\right)}(X)=
\lim\limits_{\eps\to 0}\frac{1}{2\eps}\int\limits_0^t 1_{\{x-\eps\leq X_u\leq x+\eps\}}du.
\end{equation}
In what follows we consider only  the symmetric local time (the local time).  

  Our principal result about joint density of SBM and its functionals
is the following theorem.
\begin{theorem}
\label{T1}
Let  $\left(\tau, V, L_T^{\left(0\right)}(X)\right)$ be the quantities   defined by equations (\ref{tau}), (\ref{V}) and (\ref{L})
with $X_t=W^{\left(p\right)}_t$.
Given  $X_0=0$ the joint  density of $\left(\tau, V, X_T, L_T^{\left(0\right)}(X)\right)$ is 
\begin{equation}
\label{psi-p}
\psi_{p, T}(t,v, x, l)=2a(x)h(v, lp)h(t-v, lq)h(T-t, x),\,  0\leq v\leq  t\leq T, l\geq 0, 
\end{equation} 
where $q=1-p$, 
$$a(x)=
\begin{cases}
p,& \mbox{if} \, \, x\geq 0,\\
q,& \mbox{if}\,  \, x< 0,
\end{cases}
$$
and
\begin{equation}
\label{h}
h(s, y)=\frac{|y|}{\sqrt{2\pi s^3}}e^{-\frac{y^2}{2s}},\,\, y\in \R, s\in \R_{+},
\end{equation}
is the probability density function of the  first passage time  to zero of the standard BM starting at $y$.
\end{theorem}
\begin{theorem}
\label{T4}
Let  $\left(\tau, V, L_T^{\left(0\right)}(X)\right)$ be the quantities   defined by equations (\ref{tau}), (\ref{V}) and (\ref{L})
for a solution $X_t$ of equation  (\ref{X0}).
Given  $X_0=0$ the joint  density of $\left(\tau, V, X_T, L_T^{\left(0\right)}(X)\right)$ is 
 given by the following function 
\begin{equation}
\label{phi}
\phi_{T}(t, v, x, l)
=\psi_{p, T}(t, v, x, l)e^{-\frac{m_1^2v+m^2_2(T-v)}{2}-l\left(m_1p-qm_2\right)+xm(x)},\,  0\leq v\leq  t\leq T, l\geq 0,
\end{equation}
where $\psi_{p, T}(t,v,z,l)$ is defined by equation (\ref{psi-p}).
\end{theorem}
Let us briefly comment on relationship of Theorems \ref{T1} and \ref{T4}  with some known  results.
First we  rewrite the joint density (\ref{phi}) in terms of the total occupation time.
Given $T>0$ define 
$$U=\int\limits_{0}^T 1_{\{X_t\geq 0\}}dt$$
the  total occupation time of the non-negative half-line during time period $[0, T]$
and notice  that if $X_0=0$ then 
\begin{equation}
\label{UV}
U=\begin{cases}
V+T-\tau,& \mbox{if}\,\, X_T\geq 0,\\
V,& \mbox{if}\,\, X_T<0.
\end{cases}
\end{equation}
If  $X$ is  SBM with parameter $p$ and  drift $m(x)$  then 
Theorem~\ref{T1} and  equation~(\ref{UV}) yield that 
 the joint density of  $(\tau, U, X_T, L_{T}^{\left(0\right)}(X))$
is given by the following equation 
\begin{equation}
\label{phi1}
\varphi_T(t,u, x, l)
=
\begin{cases}
2ph(u+t-T, lp)h(T-u, lq)h(T-t, x)e^{-\frac{m_1^2u+m^2_2(T-u)}{2}+xm_1-l\left(pm_1-qm_2\right)},& \\
\mbox{if} \, \, x\geq 0,\,\, l>0,\,\, \mbox{and}\,\,t\leq T,\, T-t\leq u\leq T,&\\
2qh(u, lp)h(t-u, lq)h(T-t, x)e^{-\frac{m_1^2u+m^2_2(T-u)}{2}+
xm_2-l\left(pm_1-qm_2\right)},& \\
\mbox{if}\,  \, x< 0,\,\, l>0, \,\, \mbox{and} \,\, 0\leq u\leq t\leq T.
\end{cases}
\end{equation}
If $m_1=m_2=m=const$  then we obtain the density  of the quartet  in the case of the constant drift
\begin{equation*}
\varphi_{T, m}(t,u, x, l)
=
\begin{cases}
2ph(u+t-T, lp)h(T-u, lq)h(T-t, x)e^{-\frac{m^2T}{2}+xm-lm\left(p-q\right)},& \\
\mbox{if} \, \, x\geq 0,\,\, l>0,\,\, \mbox{and}\,\,t\leq T,\, T-t\leq u\leq T,&\\
2qh(u, lp)h(t-u, lq)h(T-t, x)e^{-\frac{m^2T}{2}+
xm-lm\left(p-q\right)},& \\
\mbox{if}\,  \, x< 0,\,\, l>0, \,\, \mbox{and} \,\, 0\leq u\leq t\leq T. 
\end{cases}
\end{equation*}
Further,  setting  $m=0$ in the preceding display  and integrating  out variable $t$  we get  the joint   density 
of SBM with parameter $p$, its (total) occupation and local time (Theorem 1.2 in \cite{Appu})
\begin{align}
\nonumber
\rho(u, z, b)&=\begin{cases}
\int\limits_{0}^{T} 2ph(u+t-T, lp)h(T-u, lq)h(T-t, x)dt,& x\geq 0,\\
\int\limits_{u}^{T} 2qh(u, lp)h(T-u, lq)h(T-t, x)dt,& x<0,
\end{cases}\\
&=
\begin{cases}
2ph(T-u, bq)h(u, lp+x),&  x\geq 0,\\
2qh(u, lp)h(T-u, lq-x), & x< 0.
\end{cases}
\label{rho}
\end{align}
In a particular case $p=1/2$ density (\ref{rho}) is the trivariate density  obtained in \cite{Karatzas1} for  the standard BM.
It should be noticed that  the  local time in \cite{Karatzas1} equals to  a half of the local time defined by (\ref{L}).

\section{Application in finance}
\label{finance}

\subsection{Relationship between  LVM with discontinuity and SBM}
\label{relation}

Fix $\sigma_1>0$ and $\sigma_2>0$ and consider the following LVM
\begin{equation}
\label{S1}
dS_t=\sigma(S_t)S_tdW_t,
\end{equation}
where 
\begin{equation}
\label{sigma}
\sigma(S)=\sigma_11_{\{S\geq 1\}}+\sigma_21_{\{S<1\}}.
\end{equation}
Lemma~\ref{T0} below explains the relationship between SBM and  LVM defined by~(\ref{S1}).
 This lemma  can be regarded as a particular case of Theorem 1 in \cite{Decamps1}) (see also an argument  on~p.687 in  \cite{Decamps0})
and   is based on  application of the symmetric  Tanaka-Meyer  formula 
(e.g. see either formula (7.4) in~\cite{Karatzas3}, or  Exercise 1.25, Chapter VI in \cite{Yor},  or  formula (32)
in \cite{Lejay}).
 We  provide the proof here  for the sake of completeness and for the reader's convenience.
 \begin{lemma}
\label{T0}
Let $S_t$ be a solution of equation (\ref{S1}).
A stochastic process $X_t$ defined as follows 
\begin{equation}
\label{X1}
X_t=\frac{\log(S_t)}{\sigma(S_t)}
\end{equation}
 is a solution of the following SDE with the local time 
\begin{equation}
\label{X11}
dX_t=\mu(X_t)dt +dW_t+(p-q)dL_t^{\left(0\right)}(X),
\end{equation}
where 
\begin{equation}
\label{mu}
\mu(x)=-\frac{\sigma \left(e^{x}\right)}{2}=
\begin{cases}
 \mu _1=\left.-\sigma _1\right/2, & x\geq 0, \\
 \mu _2=\left.-\sigma _2\right/2, & x<0,\\
\end{cases}
\end{equation}
and
\begin{equation}
\label{p}
p=\frac{\sigma_2}{\sigma_1+\sigma_2},\, q=1-p=\frac{\sigma_1}{\sigma_1+\sigma_2}.
\end{equation}
In other words, $X_t$ is SBM with parameter $p=\frac{\sigma_2}{\sigma_1+\sigma_2}$ and discontinuous drift
$\mu(x)$. 
\end{lemma}
{\it Proof of Lemma  \ref{T0}.} 
First, define  $Y_t=\log(S_t)$ and notice that by usual  Ito's formula
$$dY_t=-\frac{\sigma^2(S_t)}{2}dt +\sigma(S_t)dW_t=-\frac{\sigma^2\left(e^{Y_t}\right)}{2}dt +\sigma\left(e^{Y_t}\right)dW_t.$$
In terms of process $Y_t$ we have that   $X_t=f(Y_t)$, where 
$
f(y)=\frac{y}{\sigma_1}1_{\{y\geq 0\}}+\frac{y}{\sigma_2}1_{\{y<0\}}.
$
It is easy to see that 
 $f$ is a difference of two convex functions and,  hence,  $X_t=f(Y_t)$ is a semimartingale.
Define  $f'(y)=\frac{1}{2}\left(f_l'(y)+f_r^{'}(y)\right)$, where $f_l'(y)$ and $f_r'(y)$ are 
the left and the right derivative of $f$ respectively.  
 It is easy to see  that 
 $$f^{'}(y)=\frac{1}{\sigma(y)}1_{\{y\neq 0\}}+\frac{\sigma_1+\sigma_2}{2\sigma_1\sigma_2}1_{\{y=0\}}.
 $$
 Also, the  second derivative  of $f$  (in the distribution sense) is  
$f^{''}(y)=\delta(y)\left(\frac{1}{\sigma_1}-\frac{1}{\sigma_2}\right)$, where $\delta(x)$ is the delta function.
Applying the  symmetric Tanaka-Meyer formula    to semimartingale $f(Y_t)$  we get that 
\begin{align}
\nonumber
X_t=f(Y_t)&=f(Y_0)+\int\limits_{0}^t f'(Y_u)dY_u +\int\limits_{\R}f^{''}(y)L_{t}^{\left(y\right)}(Y)dy,\\
&=f(Y_0)+\int\limits_{0}^t \left(\frac{1}{\sigma(y)}1_{\{y\neq 0\}}+\frac{\sigma_1+\sigma_2}{2\sigma_1\sigma_2}1_{\{y=0\}}
\right)dY_u
+ \frac{1}{2}\left(\frac{1}{\sigma_1}-\frac{1}{\sigma_2}\right)L_{t}^{\left(0\right)}(Y)
\label{dx}\\
 &=X_0-\int\limits_{0}^t\frac{\sigma\left(e^{X_u}\right)}{2}du +W_t+ \frac{1}{2}\left(\frac{1}{\sigma_1}-\frac{1}{\sigma_2}\right)L_{t}^{\left(0\right)}(Y),
\label{dx01}
\end{align}
where $L_{t}^{\left(0\right)}(Y)$ is the  local time of $Y_t$ at zero and where we also used that  $\int\limits_0^t 1_{\{Y_u=0\}}dY_u=0$ and   $\sigma\left(e^{Y_t}\right)=\sigma\left(e^{X_t}\right)$, in order to get equation~(\ref{dx01}) from 
equation~(\ref{dx}).

It is left  to  express $L_{t}^{\left(0\right)}(Y)$ in terms of  $L_{t}^{\left(0\right)}(X)$. 
Firstly, we apply  symmetric Tanaka-Meyer  formula  to semimartingale $X_t$  with convex  function $|x|$ and get that   
\begin{equation}
\label{f1}
|X_t|=|X_0|+\int\limits_0^t sgn(X_u)dX_u+L_t^{\left(0\right)}(X),
\end{equation}
where 
$sgn(x)=1$ if $x>0$, $sgn(x)=-1$, if $x<0$, and $sgn(0)=0$.
Secondly,  consider  $|X_t|$ as a result of applying  convex function $g(y)=|f(y)|=
\frac{y}{\sigma_1}1_{\{y\geq 0\}}-\frac{y}{\sigma_2}1_{\{y<0\}}$ to semimartingale $Y_t$. 
Let $g'$ be the arithmetic mean  of the right and the left derivatives of $g$. It is easy to see that  
$$g'(y)=\frac{1}{2}\left(g_l'(y)+g_r'(y)\right)=
sgn(y)\frac{1}{\sigma(y)}+\frac{1}{2}\left(\frac{1}{\sigma_1}-\frac{1}{\sigma_2}\right)1_{\{y=0\}}.$$
The   second generalised derivative $g''$ of $g$ is
$\left(\frac{1}{\sigma_1}+\frac{1}{\sigma_2}\right)\delta(y)$.
Applying  symmetric Tanaka-Meyer formula to $g(Y_t)$   we obtain that  
\begin{equation}
\label{f2}
|X_t|=|f(Y_t)|=|X_0|+\int\limits_{0}^t g'(Y_u)dY_u +\frac{1}{2}\left(\frac{1}{\sigma_1}+\frac{1}{\sigma_2}\right)
L_{t}^{\left(0\right)}(Y).
\end{equation}
Noticing that 
\begin{align*}
\int\limits_{0}^tsgn(X_u)dX_u
&
=\int\limits_0^t sgn(Y_u)\frac{1}{\sigma(Y_u)}\left(-\frac{\sigma^2\left(e^{Y_u}\right)}{2}du+\sigma\left(e^{Y_u}\right)dW_u\right)
\\
 &=\int\limits_0^t g'(Y_u)dY_u-g'(0)\int\limits_0^t 1_{\{Y_u=0\}}dY_u
=\int\limits_0^t g'(Y_u)dY_u
\end{align*}
and comparing right sides of equations~(\ref{f1}) and~(\ref{f2}) we  obtain  
the following identity
$$L_t^{\left(0\right)}(X)=\frac{1}{2}\left(\frac{1}{\sigma_1}+\frac{1}{\sigma_2}\right)L_{t}^{\left(0\right)}(Y),$$
which yields that 
\begin{equation}
\label{dx1}
dX_t=-\frac{\sigma\left(e^{X_t}\right)}{2} dt +dW_t+ \frac{\sigma_2-\sigma_1}{\sigma_1+\sigma_2}dL_{t}^{\left(0\right)}(X)
=\mu(X_t) dt +dW_t+(p-q)dL_{t}^{\left(0\right)}(X)\nonumber
\end{equation}
as claimed. Lemma \ref{T0} is proved.

\begin{rem}
\label{Girsanov}
{\rm 
Denote by $\Qb_T$ the probability distribution  of  SBM  with parameter $p$ and  drift~(\ref{m}) 
on the time interval $[0, T]$
 and by $\Pb_T$ the probability distribution of 
$W_t^{\left(p\right)},\, t\in [0,T]$.
By the Girsanov's  theorem  we have that
\begin{align}
\nonumber
\frac{d\Qb_T}{d\Pb_T}(X_{\cdot})&=e^{\int_{X_0}^{X_T} m(u)\, du-\frac{1}{2}\int_0^T m^2(X_t)dt
-(pm_1-qm_2)L_{T}^{\left(0\right)}(X)}\\
&=
e^{\int_{X_0}^{X_T} m(u)\, du-\frac{1}{2}\int_0^Tm^2(X_t)dt-(pm_1-qm_2)L_{T}^{\left(0\right)}(X)}\nonumber \\
&=e^{\int_{x_0}^{x} m(u)\, du-\frac{1}{2}(m^2_1w-m^2_2(T-w))-\left(pm_1-qm_2\right)l}
\label{Radon}
\end{align}
 for any trajectory 
$X_{\cdot}$ such that $X_0=x_0, X_T=x, \int_0^T1_{\{X_t\geq 0\}}dt=w,\, L_{T}^{\left(0\right)}(X)=l$.
}
\end{rem}

\subsection{Option pricing under  the two-valued  local volatility model}
\label{PricingTh}

In this section we show how the results of Section~\ref{dens} can be applied to  option pricing under the 
two-valued LVM.
We do it by example in the case of a European call option. 
 Recall first some terminology and facts of the option pricing theory. 
 A European  call option (call option) with  strike price (strike) $K$ and expiration date $T$ is a derivative whose 
 payoff is $(S_T-K)^+=\max(S_T-K, 0)$, where $S_T$ is the price of the underlying asset  at expiration. 
  A knock-in call  option with barrier $H$  is a regular call option that {\it comes into existence} only when the underlying reaches the barrier. A knock-out  call option with barrier $H$  is a regular call option that
 {\it ceases to exist}  as soon as  the underlying reaches the barrier.

Consider  the two valued LVM defined by equations  (\ref{S1}) and (\ref{sigma}) (i.e. with discontinuity at $S^*=1$).
Given value $S_0$ of the underlying at $t=0$, strike $K$ and  expiry  date $T$,
denote by  $\C=\C(S_0, K, T)$ and $\KInC=\KInC(S_0, K, T)$  the 
price of a call option and the price of a  knock-in call option with the barrier level of $1$ respectively, where both prices 
are  computed under the two-valued LVM.
Also, given the same parameters denote by 
$\KOutC=\KOutC\left(S_0, K, T, \sigma_1, 1\right)$
 the  price of a knock-out  call option  with the barrier level of $1$ computed under  the log-normal model with volatility $\sigma _1$.
 
 It is easy to see that if  $K>1$ then   
$$\C=\begin{cases}
\KInC+\KOutC,& S_0\geq 1,\\
\KInC,& S_0<1.
\end{cases}
$$
Prices of  barrier options under the log-normal model  are  known (e.g., see ch.22 in \cite{Hull}). 
Therefore, if $K>1$ then it is only left to find  $\KInC$
under the two valued LVM in order to price a call option.
A formula for the knock in call option price  $\KInC$  is given by  Theorem~\ref{T2} below.

The price of a call option with strike $K<1$ and prices of  put options can be
obtained in a similar way. Notice  that in the case of a call (put)  option with strike $K<1$ ($K>1$) 
 it seems  technically more   convenient to start with computing   the price  of a put (call) option with
the same parameters and then to  use   the put-call parity equation. 
 
\bigskip 

Let us introduce  some functions that will appear in Theorem~\ref{T2} and its proof. 
Let  
\begin{equation}
\label{pdf}
{\rm n}(x)=\frac{e^{-\frac{x^2}{2}}}{\sqrt{2\pi}},\, x\in \R,
\end{equation}
be the probability density function and
\begin{equation}
\label{cdf}
\Phi(z)=\frac{1}{\sqrt{2\pi}}\int\limits_{-\infty}^{z}e^{-\frac{y^2}{2}}dy, \, z\in \R,
\end{equation}
be  the cumulative distribution function of the standard normal distribution.
Let
\begin{equation}
\label{Norm}
\mathcal{N}\left(x,y, \rho\right)=\int\limits_{-\infty}^x\int\limits_{-\infty}^y
\frac{e^{\frac{1}{1-\rho^2} \left(-\frac{z_1^2}{2}-\rho z_1z_2+\frac{z_2^2}{2}\right)}}{2 \pi  \sqrt{1-\rho ^2}}dz_1dz_2,\, 
x,y\in \R,
\end{equation}
be the joint cumulative distribution function of the  bivariate  normal distribution   with zero means, unit variances  and correlation $\rho$.

Also denote
\begin{equation*}
\phi(S)=\frac{\log(S)}{\sigma(S)}=
\begin{cases}
 \log(S)\left/\sigma _1\right., & S\geq 1, \\
 \log(S)\left/\sigma _2\right., & S<1, 
\end{cases}
\end{equation*}
and 
\begin{equation}
\label{txb}
h(t, x, b)=\frac{x}{\sqrt{2\pi t^3}}e^{-\frac{(x+bt)^2}{2t}}, \, t\in \R_{+}, x, b\in \R,
\end{equation}
Finally, for simplicity of notation, we assume in Theorem~\ref{T2} that the risk-free interest rate is zero.
\begin{theorem}
\label{T2}
Let $S_t$ be the random process that follows equation (\ref{S1}) with function (\ref{sigma}).
Given $K>0$ and $S_0>0$ denote $k=\phi(K)$ and $x_0=\phi(S_0)$. 
Let  $\KInC=\KInC(S_0, K, T)$  be  the price of a knock-in 
European call option  with strike $K$ and  expiration date $T$ given the initial price $S_0$.
\begin{enumerate}
\item[ 1)]
If $S_0\geq 1, K>1$, then 
\begin{equation*}
\KInC=pe^{-\frac{\sigma _1x_0}{2}}\left(F_{{\rm call}}\left(\frac{\sigma _1}{2}, x_0\right)-e^{\sigma _1k}
F_{{\rm call}}\left(-\frac{\sigma _1}{2}, x_0\right)\right)
\end{equation*}
where 
\begin{align}
\label{Fcall}
F_{{\rm call}}(a, x_0)&=\int\limits_{0}^{T} F_1(T-t)F_2(a, t, x_0, 1)e^{-\frac{t\sigma^2 _1}{8}}dt,
\end{align}
and where, in turn,
\begin{align}
\label{F1}
F_1(s)&=
\frac{\sqrt{2}\left(\sigma_1e^{-\frac{\sigma^2_2 s}{8}}-\sigma_2e^{-\frac{\sigma^2 _1 s}{8}}\right)+
\sqrt{\pi  s}\sigma _1\sigma _2\left(\Phi \left(\frac{\sqrt{s} \sigma _2}{2}\right)-\Phi \left(\frac{\sqrt{s} \sigma _1}{2}\right)\right)}
{\sqrt{\pi  s}(\sigma _1-\sigma _2)},
\\
\label{F2}
F_2(a, t, x_0, \theta)&=\frac{1}{\sqrt{2 \pi } \sqrt{t}}e^{ka -\frac{\left(\left|x_0\right|+|k|\right){}^2}{2 t}}
+ae^{a\left|x_0\right|+
\frac{ta^2}{2}} \left(1-\Phi \left(\theta \frac{\left|x_0\right|+|k|}{\sqrt{t}}-a\sqrt{t}\right)\right).
\end{align}
\item[ 2)] 
If $S_0<1, K> 1$, then 
\begin{equation*}
\KInC=
2pe^{\frac{\sigma_2x_0}{2}} \left(G\left(-\frac{\sigma_1}{2}, x_0\right)-e^{k \sigma _1}
   G\left(\frac{\sigma_1}{2}, x_0\right)\right)
\end{equation*}
$$G\left(a, x_0\right)=\int\limits_0^Te^{-\frac{\sigma_1^2v}{8}-\frac{\sigma^2_2(T-v)}{8}} 
 e^{\frac{\sigma^2_1}{8}\left(u\left(\frac{p}{q}\right)^2+v\right)-|x_0|\frac{pa}{q}}  G_1\left((a, v, |x_0|, -\frac{p}{q}a\right)dv$$
where 
\begin{equation*}
G_1(a, v, y, w)=\int\limits_k^{\infty}\int\limits_0^{\infty} h(v, lp+x, a)h(T-v, lq+y, w)dldx.
\end{equation*}
In turn, $G_1$  can be expressed in terms of the standard normal distribution 
(i.e. in terms of its pdf~(\ref{pdf}) and its cdf~(\ref{cdf}))
and a  bivariate normal  cdf  (\ref{Norm}) as follows
\begin{align*}
G_1&(a, y, v, w)q\sqrt{v(T-v)}\\
&=
\frac{\n\left(\gamma X+Y\right)\n(X)}{1+\gamma^2}
-\frac{\gamma  Y}{(1+\gamma^2)^{3/2}}\n\left(\frac{Y}{\sqrt{1+\gamma^2}}\right)
\Phi\left(-\frac{(1+\gamma^2)X+\gamma  Y}{\sqrt{1+\gamma^2}}\right)\\
&-\frac{\alpha}{\sqrt{1+\gamma^2}}\,\n\left(\frac{Y}{\sqrt{1+\gamma^2}}\right)
\Phi\left(-\frac{(1+\gamma^2)X+\gamma  Y}{\sqrt{1+\gamma^2}}\right)\\
&-\beta \n(X)\Phi(-\gamma  X-Y)-
\frac{\gamma}{\sqrt{2\pi(1+\gamma^2)}}\,\n\left(\frac{Y}{\sqrt{1+\gamma^2}}\right)
 \Phi\left(-\frac{(1+\gamma^2)X+\gamma  Y}{\sqrt{1+\gamma^2}}\right)\\
&+\alpha\beta\mathcal{N}\left(-X,-\frac{Y}{\sqrt{1+\gamma^2}},-\frac{\gamma }{\sqrt{1+\gamma^2}}\right)
\end{align*}
where 
$$\alpha=w\sqrt{T-v},\, \beta=a\sqrt{v},\, \gamma=\frac{p}{q}\sqrt{\frac{T-v}{v}},$$
and 
$$X=\frac{y+(T-v)w}{\sqrt{T-v}},\, Y=\frac{qk-py-pw(T-v)+qva}{q\sqrt{v}}.$$
\end{enumerate}
\end{theorem}
Theorem \ref{T2} is proved in Section \ref{proofT2}.

\section{Proofs}
\label{proofs}

\subsection{Proofs of Theorem  \ref{T1} and Theorem \ref{T4}}

We prove Theorem \ref{T1} only. Theorem~\ref{T4} can be proved in a similar
way with straightforward modifications (see Remark~\ref{R2}). 

Given $n\in \N$  
consider   a discrete time Markov chain  $S^{\left(n\right)}_{k}\in \R,\, \, k\in \Z_{+},$ specified by the following 
transition probabilities
\begin{align*}
\P\left(\Sn_{k+1}=x+1|\Sn_{k}=x>0\right)&=\P\left(\Sn_{k+1}=x-1|\Sn_{k}=x>0\right)=\frac{1}{2},\\
\P\left(\Sn_{k+1}=x+1|\Sn_{k}=x<0\right)&=\P\left(\Sn_{k+1}=x-1|\Sn_{k}=x<0\right)=\frac{1}{2},\\
\P\left(\Sn_{k+1}=1|\Sn_{k}=0\right)&=p,\\
\P\left(\Sn_{k+1}=-1|\Sn_{k}=0\right)&=q=1-p.
\end{align*}
Define the following  stochastic process
\begin{equation}
\label{Xn}
X^{\left(n\right)}_t=\frac{1}{\sqrt{n}}\Sn_{[nt]}+\frac{nt-[nt]}{\sqrt{n}}\left(\Sn_{1+[nt]}-\Sn_{[nt]}\right),\,\, t\geq 0.
\end{equation}
Quantities  (\ref{V}), (\ref{tau}) and (\ref{L}) for  process $X^{\left(n\right)}_t$ are defined  as follows
$$\tau(\Xn)=\frac{\tau_n}{n},\,  V(\Xn)=\frac{V_n}{n},\, L(\Xn)=L_n,$$
where 
\begin{align}
\tau_n&=\max\left\{k: \Sn_k=0\right\},\label{taun}\\
V_n&=\sum\limits_{i=0}^{\tau_n}1_{\left\{\Sn_i\geq 0, \Sn_{i+1}\geq 0\right\}},\label{Vn}\\
L_n&=\sum\limits_{i=0}^{[Tn]}1_{\left\{\Sn_i=0\right\}}\label{Ln}.
\end{align}
 Theorem \ref{T1} is implied by Lemmas \ref{L0} and \ref{L1} below.
\begin{lemma}
\label{L0}
Let $X^{\left(n\right)}_t$ be the process defined by (\ref{Xn}) 
and let  $\tau_n, V_n$ and $L_n$ 
be  quantities defined by  (\ref{taun}),
(\ref{Vn}) and (\ref{Ln}).  Then 
$$\left(\frac{\tau_n}{n}, 
\frac{V_n}{n}, \frac{L_n}{\sqrt{n}}, \Xn_T\right)\to 
\left(\tau, V, L^{\left(0\right)}_T(\Wp), \Wp_T\right),$$ 
in distribution, as $n\to \infty.$
\end{lemma}
{\it Proof of Lemma \ref{L0}.}
It has been proved  in \cite{Harrison} that   
$X^{\left(n\right)}_t$ converges  in the space of continuous functions, as $n\to \infty$, 
to  SBM $\Wp_t$. This implies the claim of the lemma.

\begin{lemma}
\label{L1}
Let $X^{\left(n\right)}_t$ be the process defined by (\ref{Xn})  and 
let $\tau_n, V_n$ and $L_n$ be  quantities defined by  (\ref{taun}),
(\ref{Vn}) and (\ref{Ln}). 
Suppose that  sequences of numbers $r_n,  r_{1, n}$ and $k_n$  are   such that   
\begin{equation*}
\frac{2r_{1, n}}{n}\to x\in [0,T],\, \frac{2(r_n-r_{1,n})}{n}\to y\in [0,T], \frac{k_n}{\sqrt{n}}\to l\in \R_{+},\,\, \mbox{as}\,\, n\to \infty,
\end{equation*}
where $x+y=t\leq T$.

1) If, in addition,  $ \frac{j_n}{\sqrt{n}}\to z\geq 0$, as $n\to \infty$,
 then 
\begin{align*}
\lim\limits_{n\to \infty}\frac{n^{3}}{8}
&\P\left(V_n=2r_{1,n},\, \tau_n=2r_n, \, L_n=k_n,
\, X^{\left(n\right)}_1=j_n|X_{0}^{\left(n\right)}=0\right)\\
&=\frac{2p^2ql^2z}{[2\pi x(t-x)(T-t)]^{3/2}}e^{-\frac{z^2}{2(T-t)}-\frac{l^2}{2}\left(\frac{p^2}{x}+\frac{q^2}{t-x}\right)}.
\end{align*}
2) If, in addition, $ \frac{j_n}{\sqrt{n}}\to z< 0$, as $n\to \infty$,
 then 
\begin{align*}
\lim\limits_{n\to \infty}\frac{n^{3}}{8}&
\P\left(V_n=2r_{1,n},\, \tau_n=2r_n,  \, L^{\left(n\right)}_n=k_n,
X^{\left(n\right)}_1=j_n|X_{0}^{\left(n\right)}=0\right)\\
& =\frac{2p^2ql^2|z|}{[2\pi x(t-x)(T-t)]^{3/2}}e^{
-\frac{z^2}{2(T-t)}-\frac{l^2}{2}\left(\frac{p^2}{x}+\frac{q^2}{t-x}\right)}.
\end{align*}
\end{lemma}
Lemma \ref{L1} is proved in Section \ref{proofL1}. 

\begin{rem}
\label{R2}
{\rm  Theorem \ref{T4}  can be  proved  by  modifying appropriately  the proof of Theorem \ref{T1}.
 First of all,  the transition probabilities of  Markov chain  $S^{\left(n\right)}_{k}\in \R,\, \, k\in \Z_{+},$ should be modified 
as follows
\begin{align*}
\P\left(\Sn_{k+1}=x+1|\Sn_{k}=x>0\right)&
=\frac{1}{2}\left(1+\frac{m_1}{\sqrt{n}}\right),\\
\P\left(\Sn_{k+1}=x-1|\Sn_{k}=x>0\right)&
=\frac{1}{2}\left(1-\frac{m_1}{\sqrt{n}}\right), \\
\P\left(\Sn_{k+1}=x+1|\Sn_{k}=x<0\right)&
=\frac{1}{2}\left(1+\frac{m_2}{\sqrt{n}}\right),\\
\P\left(\Sn_{k+1}=x-1|\Sn_{k}=x<0\right)&
=\frac{1}{2}\left(1-\frac{m_2}{\sqrt{n}}\right),\\
\P\left(\Sn_{k+1}=1|\Sn_{k}=0\right)&=p,\\
\P\left(\Sn_{k+1}=-1|\Sn_{k}=0\right)&=q=1-p.
\end{align*}
Let  $\Xn$ be a stochastic  process  defined by equation~(\ref{Xn}) as before.
Convergence of $\Xn$ to SBM with drift $m$ can be proved  by a straightforward modification 
of the proof  in~\cite{Harrison} (see also~\cite{Lejay}) in the driftless case. 
Convergence implies  an analogue of Lemma~\ref{L0}. 
It is also rather straightforward to make appropriate  changes in both the statement  and the proof of Lemma \ref{L1}
in the case of non-zero drift. We skip the details. 

Alternatively, one can combine  Theorem~\ref{T1} and the Girsanov's theorem (see Remark~\ref{Girsanov}) to obtain Theorem~\ref{T4}.

}
\end{rem}

\subsection{ Proof of Lemma \ref{L1}}
\label{proofL1}
Recall that $\Xn_t$ is the process defined by (\ref{Xn}).
\begin{definition}
Given $n$ consider a discrete trajectory  $\Xn_{t_k}, k=0,1,\ldots,n$, where we denoted 
$$t_k=\frac{k}{n},\, k=0,1,\ldots,[Tn].$$
\begin{itemize}
\item A part of the trajectory    $\left(\Xn_{t_k}, \Xn_{t_{k+1}}, \ldots, \Xn_{t_{k+2d}}\right)$ such that  
$$\Xn_{t_k}=0,\,  \Xn_{t_{k+1}}>0,\,\ldots, \Xn_{t_{k+2d-1}}>0,\, \Xn_{t_{k+2d}}=0,$$
is called a positive cycle of length $2d$.

Similar, a part of the trajectory  $\left(\Xn_{t_k}, \Xn_{t_{k+1}}, \ldots, \Xn_{t_{k+2d}}\right)$ such that 
$$\Xn_{t_k}=0,\,  \Xn_{t_{k+1}}<0,\,\ldots, \Xn_{t_{k+2d-1}}<0\, \Xn_{t_{k+2d}}=0,$$
is called a negative  cycle of length $2d$.
 \item  Let $R_n$ be  the number of positive cycles  in a  trajectory $\Xn_{t_k},\, k=0,\ldots,[Tn].$
\item Given $r,  r_1, k, i\in \Z_{+}$, where $r_{1}\leq r$ and $i\leq k$, 
define  the following set of trajectories  $\Xn_{t_k},\, k=0,\ldots, [Tn],$ such that 
$$A_{r, r_1, k, i}=\left\{
\tau_n=2r,\, L_n=k,\,  R_{n}= i,\, V_n=2r_1 \right\}.$$
\end{itemize}
\end{definition}
Notice that the total number of both positive and negative cycles equals   $L_n$.
 
We prove the lemma  only if $z\geq 0$ (the case $z<0$ can be  considered  similar).
Given $j\geq 0$ denote
$$B_{n, r, j}=\left\{ X_{t_{2r+1}}^{\left(n\right)}>0, \ldots, X_{t_{[Tn]-1}}^{\left(n\right)}>0, \,  X_{T}^{\left(n\right)}=j\right\}.$$
It is easy to see that 
\begin{align*}
\P\left(V_n\right.&\left.=2r_{1,n}, \, \tau_n=2r_n, \, L_n=k_n,\, X^{\left(n\right)}_1=j_n|
X_{0}^{\left(n\right)}=0\right)\\
&=\left(\sum\limits_{i=0}^{k_n}\P(A_{r_n, r_{1,n}, k_n, i})\right)
\P\left(B_{n, r_n, j_n}|X_{t_{2r_n}}^{\left(n\right)}=0\right)
\end{align*}
The statement of the lemma is implied by two following propositions.
\begin{prop}
\label{P1}
Under assumptions of Lemma \ref{L1}
\begin{align*}
\lim\limits_{n\to \infty}
n^{2}\left(\sum\limits_{i=0}^{k_n}\P(A_{r_n, r_{1,n}, k_n, i})\right)
&=\frac{2pql^2}{\pi(x(t-x))^{3/2}}
e^{-\frac{l^2}{2}\left(\frac{p^2}{x}+\frac{q^2}{t-x}\right)}=4h(x, pl)h(t-x, lq).
\end{align*}
\end{prop}
\begin{prop}
\label{P2}
\begin{enumerate}
\item[1)] 
 Under assumptions of  Part 1) of  Lemma \ref{L1}, 
\begin{align*}
\lim\limits_{n\to \infty}n
\P\left(X_{t_{2r+1}}^{\left(n\right)}\right.&\left.>0,\,\ldots,  X_{t_{[Tn]-1}}^{\left(n\right)}>0,\, X_{T}^{\left(n\right)}=j|X_{t_{2r}}^{\left(n\right)}=0\right)\\
&=\sqrt{\frac{2}{\pi}}
\frac{2pz}{(T-t)^{3/2}}e^{-\frac{z^2}{2(T-t)}}=4ph(T-t, z).
\end{align*}
\item[2)] 
 Under assumptions of Part 2) of   Lemma \ref{L1}, 
\begin{align*}
\lim\limits_{n\to \infty}n
\P\left(X_{t_{2r+1}}^{\left(n\right)}\right.&\left.<0,\,\ldots,  X_{t_{[Tn]-1}}^{\left(n\right)}<0,\, X_{T}^{\left(n\right)}=j|X_{t_{2r}}^{\left(n\right)}=0\right)\\
&=\sqrt{\frac{2}{\pi}}
\frac{2q|z|}{(T-t)^{3/2}}e^{-\frac{z^2}{2(T-t)}}=4qh(T-t, z).
\end{align*}
\end{enumerate}
\end{prop}

\subsubsection{Proof of Proposition  \ref{P1}}

 We write $r=r_n, r_1=r_{1,n}$ and $k=k_n$ throughout  the proof.
It is easy to see that  probabilities  of a positive cycle of length $2d$ and of a negative cycle of length $2d$, where $d\geq 1$,  are equal
 to 
$2p/4^{d}$ and $2q/4^d$ respectively.  Therefore
 a probability of a  single path  from    $A_{r, r_1, k, i}$   is equal to 
\begin{equation}
\label{traj}
\frac{2^kp^{i}q^{k-i}}{2^{2r}}.
\end{equation}
 Denote by  $N_{2d, i}$ the number of paths of length $2d$, starting and  ending at the origin and 
formed by $i$ cycles  regardless of their signs. It is easy to see that the number of paths of length $2d$, starting and  ending at the origin and 
formed by $i$ cycles {\it of the same sign} is equal to $N_{2d, i}/{2^i}$. Therefore, 
 the number of trajectories forming set $A_{r, r_1, k, i}$ is equal to
\begin{equation}
\label{N}
{k\choose i} \frac{N_{2r_1, i}}{2^{i}}\frac{N_{2(r-r_1), k-i}}{2^{k-i}}.
\end{equation}
Notice that   
\begin{equation}
\label{fnk}
\frac{N_{2d, i}}{2^{2d}}=f_{2d}^{\left(i\right)},
\end{equation}
where $f_{2d}^{\left(i\right)}$ is the  probability that 
$i-$th return of  SSRW  to the origin occurs at time moment  $2d$. Summarising equations (\ref{traj}), (\ref{N}) and 
(\ref{fnk})  we get the following formula
\begin{equation*}
\P(A_{r, r_1, k, i})={k\choose i}p^{i}q^{k-i}f^{\left(i\right)}_{2r_1} f^{\left(k-i\right)}_{2(r-r_1)}.
\end{equation*}
It is known   (Section 7, ch.3, \cite{Feller})  that 
\begin{equation*}
f_{2d}^{\left(i\right)}=\frac{i}{2d-i}\frac{1}{2^{2d-i}}{2d-i\choose d}.
\end{equation*}
If $d$ is large and $i^2/(2d)$ is not very large or close to zero, then the  following approximations can be used   (equation (7.6) in Section 7, ch.3, \cite{Feller})
\begin{equation*}
f_{2d}^{\left(i\right)}\approx \sqrt{\frac{2}{\pi}}\frac{i}{(2d-i)^{3/2}}e^{-\frac{i^2}{2(2d-i)}}.
\end{equation*}
 Using this approximation it can be  shown  that  
\begin{equation}
\label{appr}
\left|\sum\limits_{i=0}^{k}{k\choose i}p^{i}q^{k-i} f^{\left(i\right)}_{2r_1} f^{\left(k-i\right)}_{2(r-r_1)}
 -\frac{2}{\pi}\sum\limits_{i=0}^{k}
\frac{ {k\choose i} p^{i}q^{k-i} i(k-i)
e^{-\frac{i^2}{2(2r_1-i)}-\frac{(k-i)^2}{2(2(r-r_1)-k+i)}}}
{(2r_1-i)^{3/2}(2(r-r_1)-k+i)^{3/2}}\right|\to 0,
\end{equation}
as $n\to \infty$.  Under assumptions of Lemma \ref{L1} the second sum in the preceding display can be replaced by the following 
one
\begin{equation}
\label{Bin}
 \frac{1}{n^2}\frac{2l^2}{\pi(xy)^{3/2}}\sum\limits_{i=0}^{k}\frac{{k\choose i}p^{i}q^{k-i}i(k-i)}{k^2}
e^{-\frac{l^2}{2}\left(\frac{i^2}{k^2x}+\frac{1}{y}\left(1-\frac{i}{k}\right)^2\right)},
\end{equation}
which, in turn,  is equal to the expectation
$\E\left(F\left(\frac{\xi_n}{k}\right)\right)$, where $\xi_n$ is a Binomial random variable with parameters $k_n$ and $p$, and 
where function $F$ is defined as follows
$$F(z)=z(1-z)e^{-\frac{l^2}{2}\left(\frac{z^2}{x}+ \frac{(1-z)^2}{y}\right)}.$$
By  Law of Large Numbers
\begin{equation}
\label{LLN}
\E\left(F\left(\frac{\xi_n}{k}\right)\right)\to F\left(p\right)=pqe^{-\frac{l^2}{2}\left(\frac{p^2}{x}+ \frac{q^2}{y}\right)},
\quad \mbox{as}\quad n\to \infty,
\end{equation}
Combining equations (\ref{appr}), (\ref{Bin}) and (\ref{LLN}) 
 we get that 
\begin{align*}
n^2 \sum\limits_{i=0}^k  \P(A_{r, r_1, k, i})
&=n^2\sum\limits_{i=0}^k{k\choose i}p^{i}q^{k-i}
 f^{\left(i\right)}_{2r_1} f^{\left(k-i\right)}_{2(r-r_1)}
\to \frac{2pql^2}{\pi(x(t-x))^{3/2}}
e^{-\frac{l^2}{2}\left(\frac{p^2}{x}+\frac{q^2}{t-x}\right)},
\end{align*}
as $n\to \infty$.

\subsubsection{Proof of Proposition \ref{P2}}

Proposition \ref{P2} is proved in \cite{Bill}, chapter 9, as   a part of derivation 
of the joint distribution of the standard BM, its last visit to the origin and the occupation time. 
 We  give the proof here  for the sake of completeness and for reader's convenience. 
For simplicity of notation and without loss of generality, we assume that $[Tn]$ is an integer, so that $t_{[Tn]}=T$.

It is easy to see that probability of a single  trajectory such that
 $$X_{t_{2r}}^{\left(n\right)}=0,\, X_{t_{2r+1}}^{\left(n\right)}>0,\,\ldots,  X_{t_{[Tn]-1}}^{\left(n\right)}>0,\, 
X_{t_{[Tn]}}^{\left(n\right)}=X_{T}^{\left(n\right)}=j>0,$$
 is equal to $p/2^{n-2r-1}$.
Therefore,
\begin{align*}
\P\left(X_{t_{2r+1}}^{\left(n\right)}\right.&\left.>0,\,\ldots,  X_{t_{Tn-1}}^{\left(n\right)}>0,\, X_{T}^{\left(n\right)}
=j|X_{t_{2r}}^{\left(n\right)}=0\right)\\
&=2p \P\left(S_{2r+1}>0,\,\ldots,  S_{Tn-1}>0,\, S_{Tn}=j|S_{2r}=0\right),
\end{align*}
where  $S_k$ is  the simple symmetric random walk  (SSRW).
 If $Tn-2r$ and $j$ have the same parity, then 
\begin{equation*}
\P\left(S_{2r+1}>0,\,\ldots,  S_{n-1}>0,\, S_{Tn}=j|S_{2r}=0\right)=\frac{j}{Tn-2r}
\P(S_{Tn-2r}=j| S_0=0).
\end{equation*}
It is easy to see that under assumptions of the lemma 
$$\frac{j}{\sqrt{Tn-2r}}\to \frac{z}{\sqrt{T-t}},$$
hence, by Local Limit  Theorem 
\begin{equation*}
\frac{\sqrt{Tn-2r}}{2}\P(S_{Tn-2r}=j| S_0=0)\to \frac{1}{\sqrt{2\pi}}e^{-\frac{z^2}{2(T-t)}}.
\end{equation*}

We conclude the proof by noticing that 
\begin{equation*}
\lim\limits_{n\to \infty}n\frac{2j}{(Tn-2r)^{3/2}}=\frac{2z}{(T-t)^{3/2}}.
\end{equation*}

\subsection{Proof of Theorem \ref{T2}}
\label{proofT2}

\paragraph{Proof of Part 1) of Theorem \ref{T2}.}
It is easy to see that if  $S_0>1$ and  $K>1$ then  we get  the following equation for the option price
 \begin{equation*}
\KInC=
\int\limits_k^{\infty}\int\limits_0^{\infty }\int\limits_{\Gamma_{T, 1}}\left(e^{\sigma _1x}-e^{\sigma _1k}\right)
e^{-\left(t_0+v+s\right) \lambda _1-u \lambda _2+\left(x-x_0\right)\mu _1}
h\left(t_0,x_0\right)
\psi_{p, T-t_0}(u+v, v, x, l)dt_0dxdldvds
\end{equation*}
where $ \lambda_i=\frac{\sigma_i^2}{8},\, i = 1, 2$, 
$t_0$ is the hitting time to zero, $v$ and $u$ are  occupation times of the positive half-line and of the negative half-line respectively 
which are observed between $t_0$ and the last visit to the origin (i.e.  $t_0+v+u$),  $s=T-(t_0+v+u)$,
$\Gamma_{T, 1}=\{(t_0, v, u, s): t_0+v+u+s=T\}$ and where $\psi_{p, T-t_0}$ is given by~(\ref{psi-p}), i.e.
$$\psi_{p, T-t_0}(u+v, v, x, l)=2ph(v, lp)h(u, lq)h(s,x),$$
since $x>0$.
  Using the  convolution property of hitting times we get that 
\begin{equation*}
\int\limits_{t_0+s=t}h\left(t_0,x_0\right) h(s,x)dt_0ds=\int\limits_{0}^t h\left(t-s,x_0\right) h(s,x)dtds
=h\left(t,\left|x_0\right|+|x|\right).
\end{equation*}
Notice  that $2ph(v,l p)h(u,l q)h\left(t,\left|x_0\right|+|x|\right)=\psi_{p, T}(v+u, v, |x_0|+|x|,  l)$ 
 and  rewrite the expression for $\KInC$  as follows
\begin{align*}
\KInC=\int\limits_k^{\infty }\int\limits_0^{\infty }\int\limits_{\Gamma_{T,2}}&
\left(e^{\sigma _1x}-e^{\sigma _1k}\right)
\psi_{p, T}(v+u,v,|x_0|+|x|, l)
e^{-(t+v)\lambda _1-u \lambda _2}e^{\mu _1\left(x-x_0\right)}dldtdvdx,
\end{align*}
where $\Gamma_{T, 2}=\{(t, v, u): t+v+u=T\}$.
Denoting 
\begin{equation*}
g(u, v)=2\int\limits_0^{\infty }h(v,l p)h(u,l q)dl=\frac{pq}{\sqrt{2 \pi } \left(p^2 u+q^2 v\right)^{3/2}}
\end{equation*}
we can rewrite 
\begin{equation}
\label{For_example}
\KInC= p\int\limits_k^{\infty }\int\limits_{\Gamma_{T, 2}}\left(e^{\sigma _1x}-e^{\sigma _1k}\right)
g(u, v) h(t, |x_0|+|x|)e^{-(t+v) \lambda _1-u \lambda _2}e^{\mu _1\left(x-x_0\right)}dtdvdx.
\end{equation}
Further, recalling  that $\sigma _1=-2\mu _1$  we arrive to the following expression for the price 
\begin{equation*}
\KInC=pe^{-\frac{\sigma _1x_0}{2}}\left(F_{{\rm call}}\left(\frac{\sigma_1}{2}, x_0\right)-e^{\sigma _1k}F_{{\rm call }}\left(-\frac{\sigma _1}{2}, x_0\right)\right),
\end{equation*}
where 
\begin{align*}
F_{{\rm call}}(\alpha, x_0)&=\int _k^{\infty }\int\limits_{t+v+u=T} g(u, v)
e^{-v \lambda _1-u \lambda _2}h\left(t,\left|x_0\right|+|x|\right)e^{\alpha x}e^{-t \lambda _1}dtdvdx.
\end{align*}
Integrating with respect to variables $u,v$ provided that  $u+v=T-t=s$ is fixed we obtain function
\begin{align*}
F_1(s)&=\int\limits_{v+u=s}g(u, v)e^{-v \lambda _1-u \lambda _2}dv\\
&=\frac{\sigma _1\sigma _2 }{\sigma_1-\sigma _2}\left[\frac{2}{\sqrt{2 \pi  s}}
\left(\frac{e^{-\frac{1}{8} s \sigma _2^2}}{\sigma _2} -\frac{e^{-\frac{1}{8} s \sigma_1^2}}{\sigma _1}\right)+
\left(\Phi \left(\frac{\sqrt{s} \sigma _2}{2}\right)-\Phi \left(\frac{\sqrt{s} \sigma _1}{2}\right)\right)\right],
\end{align*}
defined earlier by equation~(\ref{F1}). Integrating out variable $x$   we get
\begin{align*}
\int _k^{\infty }h\left(t,\left|x_0\right|+|x|\right)e^{\alpha  x}dx&
=\frac{1}{\sqrt{2 \pi } \sqrt{t}}e^{k \alpha -\frac{\left(\left|x_0\right|+|k|\right){}^2}{2t}}\\
&+
\alpha  e^{\frac{t \alpha ^2}{2}-\alpha  \left|x_0\right|}\left(1-\Phi \left(\frac{|x_0|+|k|-t \alpha }{\sqrt{t}}\right)\right)\\
&=F_2(\alpha, t, x_0, 1),
\end{align*}
where function  $F_2(\alpha, t, x_0, \theta)$  is   defined by~(\ref{F2}). 
Finally, we rewrite  $F_{{\rm call}}$ in terms of $F_1$ and $F_2$ 
$$
F_{{\rm call}}(\alpha, x_0)=\int\limits_{0}^{T} F_1(T-t)F_2(\alpha, t, x_0,  1)e^{-\frac{t\sigma_1^2}{8}}dt,
$$
as claimed in~(\ref{Fcall}).

\paragraph{ Proof of Part 2) of Theorem \ref{T2}.} 
 If $S_0<1$ and $K>1$, then
 $x_0=\phi(S_0)=\frac{\log(S_0)}{\sigma _2}<0, k=\phi(K)=\frac{\log(K)}{\sigma _1}>0$,
 and we get, using notation introduced in the proof of Part 1), that 
\begin{equation*}
\KInC=\int\limits_k^{\infty }\int\limits_0^{\infty }\int\limits_{\Gamma_{T, 1}}
\left(e^{\sigma _1x}-e^{\sigma _1k}\right)h\left(t_0,x_0\right)
\psi_{p, T-t_0}(v+u, v, x, l)
e^{- \lambda_1(v+s)- \lambda_2(t_0+u)+\mu _1x-\mu_2x_0}dt_0dxdldvds,
\end{equation*}
where, as before, $h(t_0, x_0)\psi_{p, T-t_0}(u+v, v, x, l)=2ph(t_0, x_0)h(v, lp)h(u, lq)h(s,x).$
We  use again   the  convolution property of hitting times  as in Part 1) but integrate now products
$h(v, lp)h(s, x)$  and $h(t_0, x_0)h(u, lq)$ given constraints   $v+s=const$ and $t_0+u=const$ respectively.
It leads to the following expression for the price 
\begin{align*}
\KInC&=2p\int\limits_0^T\int\limits_k^{\infty }\int\limits_0^{\infty }
\left(e^{\sigma _1x}-e^{\sigma _1k}\right)h(v,lp+x)h(u,lq+|x_0|)
e^{-\lambda _1v-\lambda _2u+\mu _1x-\mu _2x_0}dldxdv\\
&=2pe^{\frac{\sigma_2x_0}{2}}\int\limits_0^Te^{-\frac{\sigma^2_1v}{8}-\frac{\sigma^2_2u}{8}}\int\limits_k^{\infty }\int\limits_0^{\infty }
\left(e^{\sigma _1x}-e^{\sigma _1k}\right)h(v,lp+x)h(u, lq+|x_0|)e^{-\frac{\sigma_1x}{2}}dldxdv\\
&=2pe^{\frac{\sigma_2x_0}{2}}\int\limits_0^Te^{-\frac{\sigma^2_1v}{8}-\frac{\sigma^2_2u}{8}}\int\limits_k^{\infty }\int\limits_0^{\infty }
h(v,lp+x)h(u,lq+|x_0|)e^{\frac{\sigma _1x}{2}}dldxdv\\
&-2pe^{\frac{\sigma_2x_0}{2}+\sigma _1k}\int\limits_0^Te^{-\frac{\sigma^2_1v}{8}-\frac{\sigma^2_2u}{8}}
\int\limits_k^{\infty }\int\limits_0^{\infty }
h(v,lp+x)h(u,lq+|x_0|)e^{-\frac{\sigma_1x}{2}}dldxdv
\end{align*}
where   $u=T-v$ and  $\mu_i=-\sigma_i/2$ and $\lambda_i=\sigma^2_i/8$.
Rewrite 
\begin{equation}
\label{C2}
 \KInC=
2pe^{\frac{\sigma_2x_0}{2}}\int\limits_{0}^Te^{-\frac{\sigma^2_1v}{8}-\frac{\sigma^2_2(T-v)}{8}}\left({\rm I}\left(-\frac{\sigma_1}{2},  |x_0|, v\right)
-e^{\sigma _1k}{\rm I}\left(\frac{\sigma_1}{2}, |x_0|, v\right)\right)dv,
\end{equation}
where  
\begin{equation}
\label{Int}
{\rm I}(a, y, v)=\int\limits_k^{\infty }\int\limits _{0}^{\infty}h(v, lp+x)h(u, lq+y)e^{-ax}dldx,\quad y\geq 0.
\end{equation}
Notice that  
\begin{align*}
h(v, lp+x)h(u, lq+y)e^{-ax}&=h(v, lp+x, a)e^{\frac{va^2}{2}+alp}h(u, lq+y)\\
&=e^{\frac{va^2}{2}+\frac{u}{2}\left(\frac{ap}{q}\right)^2-ay\frac{p}{q}}h(v, lp+x, a)h\left(u, lq+y,-apq^{-1}\right)
\end{align*}
where 
$h(t, x, b)$ is defined by~(\ref{txb}), so we can rewrite 
$${\rm I}(a, y, v)=e^{\frac{va^2}{2}+\frac{T-v}{2}\left(\frac{ap}{q}\right)^2-ay\frac{p}{q}}
G_1\left(a, v, y,  -apq^{-1}\right),$$
where 
\begin{equation}
\label{G1}
G_1(a,  v, y, w)=\int\limits_k^{\infty }\int\limits _{0}^{\infty}h(v, lp+x, a)h(u, lq+y,w)dldx.
\end{equation}
It is left  to show that $G_1(a, y, v, w)$ can be expressed in terms of the standard normal distribution and a
bivariate normal distribution.
Noticing  that 
$$h(v, lp+x, a)h(u, lq+y,w)=\frac{(lp+x)(lq+y)}{2\pi (uv)^{3/2}}e^{-\frac{(lp+x+va)^2}{2v}-\frac{(lq+y+uw)^2}{2u}}$$
and changing  variables
$z_1=\frac{lp+x+av}{\sqrt{v}}$ and $z_2=\frac{lq+y+uw}{\sqrt{u}}$ 
we  can rewrite the expression for $G_1$ as follows 
$$G_1(a, y, v, w)=\int\limits_{D}\frac{e^{-\frac{w_1^2}{2}-\frac{w_2^2}{2}}}{2\pi}
\frac{(z_1-a\sqrt{v})(z_2-w\sqrt{u})}{q\sqrt{uv}}dz_1dz_2$$
where 
$D=\{(z_1, z_2)\in \R^2:z_2\sqrt{u}>y+uw,\, -z_2p\sqrt{u}+q\sqrt{v}w_2>qk-py+qva-puw\}.$
Denote
$$\alpha=w\sqrt{u},\, \beta=a\sqrt{v},\, X=\frac{y+uw}{\sqrt{u}},\, 
Y=\frac{qk-py-puw+qva}{q\sqrt{v}},\, \gamma=\frac{p}{q}\sqrt{\frac{u}{v}},$$
and 
$\Gamma=\{(z_1, z_2):  z_1>Y+\gamma z_2,\,  z_2>X\}.$
In these notation
$$G_1(a, y, v, w)=\int\limits_{\Gamma}\frac{e^{-\frac{z_1^2}{2}-\frac{z_2^2}{2}}}{2\pi}
\frac{(z_1-\beta)(z_2-\alpha)}{q\sqrt{uv}}dz_1dz_2
=\frac{1}{q\sqrt{uv}}{\rm J}(a, y, v, w),
$$
where ${\rm J}(a, y, v, w)=\int_{\Gamma}
\n(z_1)\n(z_2)
(z_1-\beta)(z_2-\alpha)dz_1dz_2,$
and function $\n$ is  defined  by~(\ref{pdf}).
Notice that 
\begin{align*}
{\rm J}(a, y, v, w)&=\int\limits_{\Gamma}z_1z_2\n(z_1)\n(z_2)
dz_1dz_2
-\alpha\int\limits_{\Gamma}z_1\n(z_1)\n(z_2)
dz_1dz_2\\
&-\beta\int\limits_{\Gamma}z_2\n(z_1)\n(z_2)
dz_1dz_2+\alpha\beta\int\limits_{\Gamma}\n(z_1)\n(z_2)dz_1dz_2\\
&:={\rm J}_{1}+{\rm J}_{2}+ {\rm J}_{3}+{\rm J}_{4}.
\end{align*}
It can be shown (we skip intermediate computational details) that  
\begin{align*}
{\rm J}_{1}
&=\frac{\n\left(\gamma X+Y\right)\n(X)}{1+\gamma^2}
-\frac{\gamma  Y}{(1+\gamma^2)^{3/2}}\n\left(\frac{Y}{\sqrt{1+\gamma^2}}\right)
\Phi\left(-\frac{(1+\gamma^2)X+\gamma  Y}{\sqrt{1+\gamma^2}}\right)\\
 {\rm J}_{2}
&=-\frac{\alpha}{\sqrt{1+\gamma^2}}\,\n\left(\frac{Y}{\sqrt{1+\gamma^2}}\right)
\Phi\left(-\frac{X(1+\gamma^2)+\gamma  Y}{\sqrt{1+\gamma^2}}\right)\\
{\rm J}_{3}
&=-\beta \n(X)\Phi(-\gamma  X-Y)-
\frac{\gamma}{\sqrt{2\pi(1+\gamma^2)}}\,\n\left(\frac{Y}{\sqrt{1+\gamma^2}}\right)
 \Phi\left(-\frac{(1+\gamma^2)X+\gamma  Y}{\sqrt{1+\gamma^2}}\right)\\
{\rm J}_{4}&=\alpha\beta\mathcal{N}\left(-X,-\frac{Y}{\sqrt{1+\gamma^2}},-\frac{\gamma }{\sqrt{1+\gamma^2}}\right)
\end{align*}
This finishes the proof of the second part of the theorem.

\section{A Black-Scholes approximation}
\label{simple}

In this section we  derive   a  surprisingly  simple and accurate approximation for the option price which is based on the  Black-Scholes (BS) formula.
 We use the same notation as in Sections~\ref{PricingTh} and~\ref{proofT2}.
Without loss of generality assume  that $S_0=1$ ($x_0=0$) and $K>1$ ($k>0$).  In this case $\C=\KInC$ 
and   equation~(\ref{For_example}) becomes
\begin{equation}
\label{Cin}
\C= p\int\limits_k^{\infty }\int\limits_{t+u+v=T}\left(e^{\sigma _1x}-e^{\sigma _1k}\right)
g(u, v) h(t, x)e^{-(t+v) \lambda _1-u \lambda _2}e^{\mu _1x}dtdvdx.
\end{equation}
The approximation is motivated by the following idea.  
Since $k>0$ we "should be mostly interested" in those trajectories of  $X_t$ that  spend the "most of their  lifetime" 
in region $X_t>0$, where $\sigma=\sigma_1$. Therefore, let us first replace 
 function $e^{-(t+v) \lambda _1-u \lambda _2}$  in~(\ref{Cin})  
 by $e^{-\lambda _1T}$.  Secondly,  integrating  out variables $v$ and $u=T-t-v$ gives 
$$\int\limits_{0}^{T-t}g(u, v)dv=\int\limits_0^{T-t}\frac{pq}{\sqrt{2\pi}(p^2(T-t-v)+q^2v)^{3/2}}dv
=\frac{2}{\sqrt{2\pi(T-t)}}=2{\rm p}(0, T-t),$$
 where ${\rm p}(y, T-t)$ is the transition density of the standard BM at time $T-t$ so that the result of  integration 
does not  depend on $p$ and $q$. 
Thus, we arrive, after 
expressing both $\lambda_1$ and $\mu_1$ in terms of $\sigma_1$,   to  
  the following approximation for the option price 
\begin{align}
\nonumber
\C&\approx 2p\int\limits_k^{\infty }\left(e^{\sigma _1x}-e^{\sigma _1k}\right)e^{- \frac{\sigma^2 _1T}{8}}e^{-\frac{\sigma _1x}{2}}
\left(\int\limits_{0}^Th(t, x){\rm p}(0, T-t)dt\right)dx\\
\nonumber 
&= 2p\int\limits_k^{\infty }\left(e^{\sigma _1x}-e^{\sigma _1k}\right)
e^{- \frac{\sigma^2 _1T}{8}-\frac{\sigma _1x}{2}}{\rm p}(x, T)dx\\
&=\frac{2\sigma_2}{\sigma_1+\sigma_2}
\BSCall(\sigma_1) \label{BScall}
\end{align}
where  $\BSCall(\sigma_1)$ is the BS  price of the option  under the log-normal model with volatility $\sigma_1$.
It is obvious that if  we set $\sigma_1=\sigma_2$ in both sides of the preceding display, 
 then the approximation would become the  BS formula for the call option price with volatility $\sigma_1$.

Using the same argument we can obtain similar approximation for the put option price. Namely,  if $S_0=1$, then
the price of a put option with strike $K<1$ can be approximated as follows 
\begin{equation}
\label{BSput}
\Put \approx 2q\BSPut(\sigma_2)=\frac{2\sigma_1}{\sigma_1+\sigma_2}\BSPut(\sigma_2),
\end{equation}
where  $\BSPut(\sigma_2)$  is the BS  price of the put  option  with volatility $\sigma_2$.
 Similar to the case of the call option, the BS approximation provides 
 either an upper bound (if $\sigma_1>\sigma_2$) or a lower bound (if $\sigma_1<\sigma_2$).

\begin{figure}[htbp]
\centering
\begin{tabular}{cc}
   \includegraphics[width=2.8in, height=2.0in]{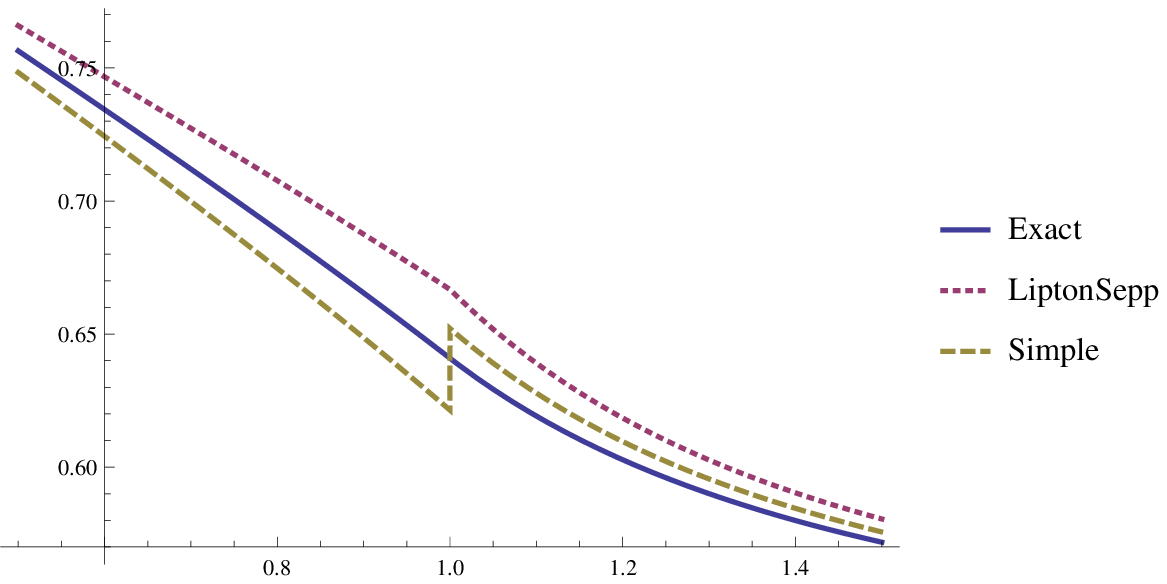}&
     \includegraphics[width=2.8in, height=2.0in]{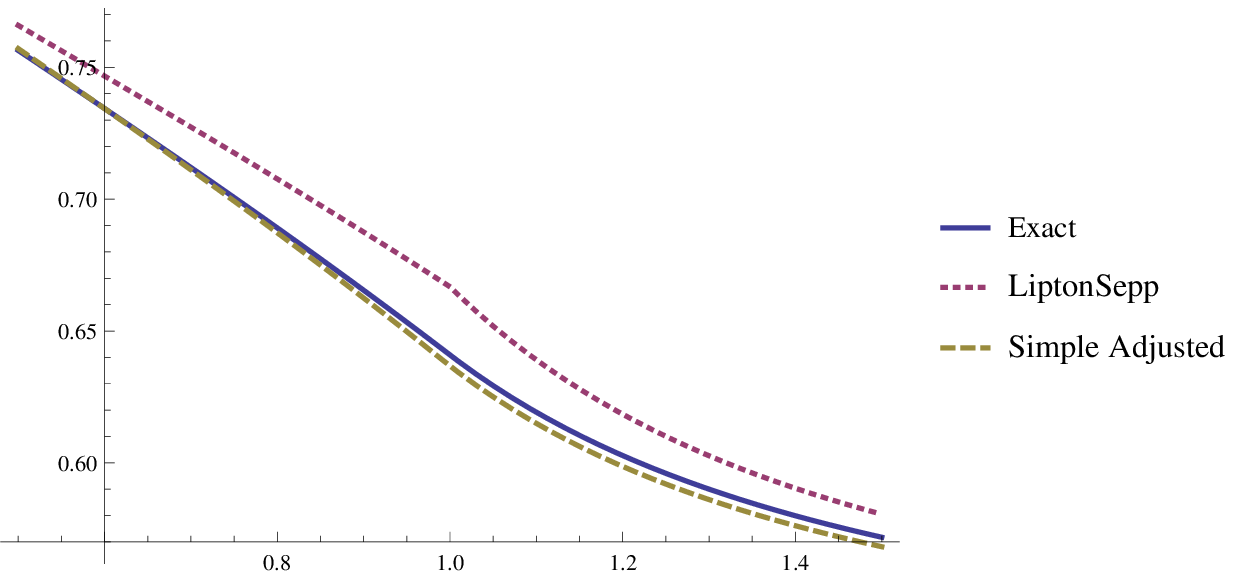} \\
\end{tabular}
\caption{{\footnotesize   Implied volatility curves, $\sigma_1=0.5,\, \sigma_2=0.9, T=2, S_0=1$. In both figures:
the solid line corresponds to the two-valued LVM and  the dashed upper curve corresponds to Lipton-Sepp's approximation.
Implied volatility calculated by using BS approximation: without adjustment on the left and with adjustment on the right.
}}
\label{Fig1}
\end{figure}

A discontinuous (at $K=1$) curve    shown on the left-hand side of Figure~\ref{Fig1} 
 is the implied volatility curve calculated by using the approximation.
In this calculation  call prices have been used, if $K>1$, and put prices have been used, if $K<1$.
A solid curve in the middle of the left side Figure~\ref{Fig1}  is the implied volatility curve calculated by using the exact formula provided by Theorem~\ref{T2}.
  It is easy to see that if $\sigma_1<\sigma_2$ then the BS approximation provides an upper (lower) bound of the price in the case of call (put) options,
and, vice versa, if $\sigma_1>\sigma_2$ then the  approximation  provides a  lower (upper)  bound for call  (put) prices.
In this example    $\sigma_1=0.5<\sigma_2=0.9$, therefore  the approximate curve is below the exact curve, if $K<1$, and above it, 
if $K>1$, as expected.
The  upper dashed curve 
 is the implied volatility curve calculated  by using 
an approximation proposed  in~\cite{Sepp} for calibration of a LVM with a piecewise volatility (tiled LVM). The latter includes
 the  two-valued LVM  as a particular case.

The  BS approximation can be improved. Indeed, 
 recall that we must have  the put-call parity  $\C-\Put=K-S_0$, which becomes $\C=\Put$, if $K=S_0$.
The put-call parity does not hold for the approximate prices and we 
adjust them so that  the put-call parity holds at $K=1$.
Namely,  define the following adjusting factors
$$A_{{\rm cl}}=\frac{p\BSCall(\sigma_1)+q\BSPut(\sigma_2)}{2p\BSCall(\sigma_1)},\,\,\, A_{{\rm pt}}=\frac{p\BSCall(\sigma_1)+q\BSPut(\sigma_2)}{2q\BSPut(\sigma_2)},$$
and redefine the approximate prices as $\widetilde{\BSCall(\sigma_1)}=A_{{\rm cl}}\BSCall(\sigma_1)$ 
and $\widetilde{\BSPut(\sigma_2)}=A_{{\rm pt}}\BSPut(\sigma_2).$ 
By construction, the put-call parity now holds for adjusted prices at $K=1$. This adjustment  smooths  the approximate  implied volatility curve
which becomes continuous everywhere.
The result of adjustment is  shown on  the right-hand  side of Figure~\ref{Fig1}, where both  the solid line and  the upper dashed line  are as before, and  the new dashed  curve is calculated by using adjusted prices.  
It is quite visible that the adjustment improves the  approximation.

Finally,   numerical tests showed that  accuracy of the approximation  improves as the time to  expiration 
 becomes smaller, which agrees with  intuition.

\section{A note on a displaced diffusion model with discontinuity }
\label{displace}

Our results on the joint distribution of SBM and its functionals can be also   applied to  derivative 
pricing in  the following displaced   model 
\begin{equation}
\label{Loc1}
dS_t=\left(\sigma_1\left(S_t-\alpha _1\right)1_{\{S_t\geq S^*\}}+\sigma_2\left(S_t-\alpha _2\right)1_{\{S_t<S^*\}}\right)dW_t,
\end{equation}
where $\sigma_1\neq \sigma_2,\,  \alpha_i\in \R,\, i=1,2$ and $S^*>0$.
Model~(\ref{Loc1})   is a particular case of the following model considered in~\cite{Decamps1}
\begin{equation*}
dS_t=\left(\sigma_1\left(S_t-\alpha _1\right)^{\beta_1}1_{\{S_t\geq S^*\}}+\sigma_2\left(S_t-\alpha _2\right)^{\beta_2}
1_{\{S_t<S^*\}}\right)dW_t.
\end{equation*}
where, in addition, $\beta_i\geq 0,\, i=1,2$.
In~\cite{Decamps1} they derived certain semi-analytical  expressions  for  the transition density of the underlying process.
The technique of~\cite{Decamps1} is an adaptation of a technique that was used in~\cite{Gorovoi}.
In turn, the technique of~\cite{Gorovoi}  is based on a well known  observation  (e.g. \cite{Gikhman}) that 
the transition density satisfies a partial differential equation and can be constructed by means of an eigenfunction expansion in the corresponding 
Sturm-Liouville problem.  In general, these  eigenfunction expansions  for the transition densities  are difficult to handle  analytically and 
an  approximation is  required. 
It should be noticed that  in~\cite{Decamps1} an  analytical expression for the  transition density was obtained in 
 a particular case where $\sigma_1=\sigma_2,\, \alpha_1\neq \alpha_2$ so that 
dependence of the joint density  on the occupation time
becomes trivial  (e.g., see equation~(\ref{phi}) or~(\ref{Radon}), where $m_1=m_2$).

Notice also that if  $\sigma_1=\sigma_2=\sigma, \,  \beta _1=\beta _2=1$ and $\alpha _1=\alpha _2=a$ then it is a classical case of a displaced log-normal 
model. The latter is just $S_t=Z_t-a$, where $Z_t$ is the  log-normal process, 
and it can be written in the local volatility form, namely,
$dS_t=\sigma(1-a/S_t)S_tdW_t$. 
 The displaced diffusion is a very useful tool for approximating  more complicated stochastic processes  in finance. 
The main reason is that this model is a first-order approximation of any LVM 
(see Remark 7.2.14 in~\cite{Piterbarg}   and other examples therein). 
A known  problem with a  displaced model of any sort is  that theoretically the underlying process can take negative values 
(e.g. when $\alpha_i>0$). 
This problem can be dealt with by imposing some constraints. For instance, instead of the classic displaced log-normal model one can consider 
model~(\ref{Loc1}) with \(\alpha _2=0\). This means that the  volatility is  a  hyperbolic function 
 above level $S^*$ and a constant one  below  level $S^*$ and, hence,  is prevented to take large values as the process approaches $0$. 
It is rather straightforward to apply our results to the displaced log-normal  model with such constraints. 
Let us take, for example,  model~(\ref{Loc1}), where  $S^*=1,  \alpha_1<1$ and $\alpha_2=0$, and consider briefly 
 the case when the process starts at $S_0<1$.  Given $\sigma_1, \sigma_2, \alpha_1$ and strike $K>1$ define
$$p=\frac{\sigma _2}{\sigma _2+\sigma _1\left(1-\alpha _1\right)},\,q=1-p, \, k=\frac{1}{\sigma_1}\log\left(\frac{K-\alpha _1}{1-\alpha _1}\right),\, 
x_0=\frac{\log\left(S_0\right)}{\sigma_2},
b =\frac{q\sigma _2-p\sigma _1}{2}
.$$
Then the price of a knock-in European call option with strike $K$ and expiration date $T$ is given by the following integral
\begin{align*}
\KInC=2p\left(1-\alpha _1\right)\int\limits_{k}^{\infty}\int\limits_{0}^{\infty}\int\limits_{\Gamma_{T, 1}}
&\left(e^{\sigma _1 x}-e^{\sigma_1k}\right)e^{-l\beta - \lambda _1 (s+v)-x_0\mu_2-\lambda_2 \left(t_0+u\right)+\mu _1 x}
R(u, v, x, l,  t_0)dxdldvdt_0
\end{align*}
where $R(u, v, x, l,  t_0)=h\left(t_0, x_0\right)\psi_{p, T-t_0}(u+v, v, x, l)$ and we used notation introduced in  the proof of Part 1)
of Theorem~\ref{T2}.
Using the same argument as in the proof of the theorem  one can show that computation of  the above integral 
can be reduced to computation of the following integral
$$\tilde{{\rm I}}(b, a, v, y)=\int\limits_k^{\infty }\int\limits_0^{\infty }e^{-ax-b l} h(v,l p+x) h(u,l q+y)dldx.$$
In turn, one can express, by modifying appropriately the argument applied to integral~(\ref{Int}),
the integral in the preceding display  
  in terms of both a univariate and a bivariate normal distribution   as follows
\begin{align*}
\tilde{{\rm I}}(b, a, v, y)&=\frac{e^{\frac{\nu^2 u}{2}+\frac{a^2 v}{2}+\nu y}}{q \sqrt{u v}}
\left(\frac{e^{-\frac{X^2+(Y+\gamma X)^2}{2}}}{2 \pi  \left(1+\gamma^2\right)}
-\frac{B e^{-\frac{X^2}{2}} \Phi(-Y-\gamma X)}{\sqrt{2\pi }}\right.\\
&+\frac{1}{\sqrt{2\pi(1+\gamma^2)}}\left(-A+B \gamma -\frac{\gamma  Y}{1+\gamma^2}\right)
e^{-\frac{Y^2}{2 \left(1+\gamma^2\right)}}\Phi\left(-\frac{(1+\gamma^2)X+\gamma Y)}{\sqrt{1+\gamma^2}}\right)\\
&\left.
+AB\mathcal{N}\left(-X,-\frac{Y}{\sqrt{1+\gamma^2}},-\frac{\gamma }{\sqrt{1+\gamma^2}}\right)\right)
\end{align*}
where $\nu=-\frac{a p-b }{q}$, $A=\nu\sqrt{u}$, $B=a\sqrt{v}$ , $Y=\frac{q (k+a v)-p(\nu  u+y)}{q \sqrt{v}},$
$X=\frac{\nu  u+y}{\sqrt{u}}$ and $\gamma=\frac{p}{q}\sqrt{\frac{u}{v}}$.


\begin{thebibliography}{1}
\bibitem{Piterbarg}   Andersen L. and Piterbarg V. (2010).  Interest Rate Modeling. Atlantic Financial Press.
\bibitem{Appu1}   Appuhamillage, T.,  Bokil, V.,  Thomann, E.,  Waymire, E.,  Wood, B. (2011):
 Corrections: Occupation and Local Times for Skew Brownian Motion with Applications to Dispersion Across an Interface. {\it Annals of  Applied
 Probability}, {\bf 21}, N5, pp.~2050--2051.
\bibitem{Appu}   Appuhamillage, T.,  Bokil, V.,  Thomann, E.,  Waymire, E.,  Wood, B. (2011). 
 Occupation and local times for skew Brownian motion with applications to dispersion across an interface. {\it Annals of  Applied
 Probability}, {\bf 21}, N1, pp.~183--214.
\bibitem{Bill} Billingsley P., (1968). {\it Convergence of Probability Measures.}  John Wiley$\&$Sons, Inc.
\bibitem{Decamps0} Decamps, M., De Schepper, A. and Goovaerts, M. (2004). Applications of $\delta-$function perturbation
to the  pricing of derivative securities. {\it Physica A}, {\bf 342}, pp.~677--692.
\bibitem{Decamps1} Decamps, M.,  Goovaerts, M., and Schoutens, W. (2006). Self Exciting Threshold Interest Rates Models.
{\it International Journal of Theoretical and Applied Finance}, {\bf  9}, N7, pp.~1093-1122. 
\bibitem{Decamps2} Decamps, M.,  Goovaerts, M., and Schoutens, W. (2006). Asymmetric skew Bessel processes and their
applications to finance. {\it Journal of Computational and Applied Mathematics}, {\bf  186},  pp.~130--147.
\bibitem{Dupire}  Dupire, B.  (1994). Pricing with a Smile. {\it Risk}, 7(1), pp.~18--20. 
\bibitem{Feller} Feller, W. (1968).  An Introduction to Probability Theory and its Applications. Volume 1, 3rd Edition. John Wiley$\&$Sons, Inc.
\bibitem{Gikhman} Gikhman, I., and  Skorohod, A. (1972). Stochastic differential equations.
\bibitem{Gorovoi} Gorovoi, V. and  Linetsky, V. (2004). Black’s model of interest rates as
options, eigenfunction expansions and Japanese interest rate. {\it Mathematical finance}, {\bf 14}, N1,  pp.~49--79.
\bibitem{Guyon} Guyon, J. and Henry-Labordere, P. (2011). From spot volatilities to implied volatilities. {\it Risk magazine}, {\bf 5}.
\bibitem{Gyongy}  Gyongy, I.  (1986): Mimicking the One-Dimensional Marginal Distributions of Processes Having an Ito Differential.
{\it Probability Theory and Related Fields}, {\bf 71}, pp.~501--516.
\bibitem{Harrison} Harrison, J. and Shepp, L. (1981). On skew Brownian motion. {\it Annals of  Probability}, {\bf 9}, N2, pp.~309--313.
\bibitem{Hull} Hull, J. (2009). Options, Futures and Other Derivatives. Pearson/Prentice Hall.
\bibitem{Karatzas1} Karatzas, I., and Shreve, (1984). Trivariate density of Brownian motion, its local and occupation times, with application to stochastic
control. {\it Annals of Applied Probability}, {\bf 12}, N3, pp.~819--828.
\bibitem{Karatzas3} Karatzas, I., and Shreve, S. (1991). Brownian motion and Stochastic Calculus. Springer-Verlag, New York, 2nd Edition. 
\bibitem{Ito} Ito, K.  and  McKean, H. (1963). Brownian motions on a half line. 
{\it Illinois  Journal of  Mathematics,} {\bf 7}, pp.~181--231.  
\bibitem{LeGall}  Le Gall, J.-F. (1985). One-dimensional stochastic differential equations involving the local times of the  unknown process.
In {\it Stochastic Analysis and Applications.  Lecture Notes in Math.}, {\bf 1095}. Springer-verlag, pp.~51-82. 
\bibitem{Lejay} Lejay, A. (2006). On the constructions of the
skew Brownian motion. {\it Probability Surveys}, {\bf 3}, pp.~413--466.
\bibitem{Sepp} Lipton, A. and Sepp, A. (2011). Filling the Gaps. {\it Risk},  October, pp.~86--91.
\bibitem{Lulko} Lyulko, Ya. (2012). On the Distribution of Time Spent by a Markov Chain at Different
 Levels Until Achieving a Fixed State. {\it Theory of Probability and  its Applications}, {\bf 56}, N1, pp.~140--149. 
\bibitem{Yor} Revuz, D. and Yor, M. (1998). {\it Continuous martingales and Brownian motion.} Springer-Verlag,  3rd Edition. 
\bibitem{Rossello} Rossello, D. (2012). Arbitrage in Skew Brownian Motion Models. 
{\it Insurance Mathematics and Economics}, {\bf 50}, N1, pp.~50--56.
\bibitem{Tak} Takacs, L.   (1995). On the local time of Brownian motion. {\it Annals of  Applied  Probability}, {\bf 5}, N3, pp.~741--756.
\end{thebibliography}
\end{document}